\documentclass[a4paper,10pt]{amsart}

\usepackage[utf8x]{inputenc}
\usepackage{amssymb}
\usepackage{amsthm}
\usepackage{amsmath}
\usepackage{amsfonts}
\usepackage{bbm}
\usepackage[T1]{fontenc}
\usepackage{graphicx}
\usepackage[]{overpic}
\usepackage{import}
\usepackage[dvips]{epsfig}
\usepackage[english]{babel}
\usepackage{array}
\graphicspath{{images/}}
\usepackage[a4paper]{geometry}
\usepackage{color}
\usepackage{setspace}
\usepackage{fancybox}
\usepackage{hyperref}
\usepackage{pinlabel} 
\usepackage{xcolor}
\definecolor{green}{RGB}{0,150,0}

\newcommand{\Z}{\mathbb Z} 

\newcommand{\R}{\mathbb R}

\newcommand{\La}{\Lambda}

\newtheorem{teo}{Theorem}

\theoremstyle{definition}
\newtheorem{defi}{Definition}

\newtheorem{rem}{Remark}

\begin{document}
	\title{Doubly slice knots and obstruction to Lagrangian concordance.}
\author{Baptiste Chantraine}
\author{No\'emie Legout}
\address{ Nantes Université, CNRS, Laboratoire de Mathématiques Jean Leray, LMJL,
UMR 6629, F-44000 Nantes, France.}
\email{baptiste.chantraine@univ-nantes.fr}

\address{Uppsala University, Department of Mathematics, Box 480, 751 06 Uppsala, Sweden.}
\email{noemie.legout@math.uu.se}
	\date{}
	\maketitle
	\begin{abstract}
          In this short note we observe that a result of Eliashberg and Polterovitch allows to use the doubly slice genus as an obstruction for a Legendrian knot to be a slice of a concordance from the trivial Legendrian knot with maximal Thurston-Bennequin invariant to itself. This allows to obstruct concordances from the Pretzel knot $P(3,-3,-m)$ when $m\geq 4$ to the unknot. Those examples are of interest because the Legendrian contact homology algebra cannot be used to obstruct such a concordance.
	\end{abstract}
	
\section{Introduction and results}

\begin{defi}
	Let $\La^-,\La^+$ be two Legendrian knots in $S^3$. A Lagrangian concordance from $\La^-$ to $\La^+$ is a Lagrangian submanifold $\Sigma\in\R\times S^3$ diffeomorphic to a cylinder and such that for some $T>0$,
	\begin{itemize}
		\item $\big((-\infty,-T)\times S^3\big)\cap\Sigma=(-\infty,-T)\times\La^-$ and
		\item $\big((T,+\infty)\times S^3\big)\cap\Sigma=(T,+\infty)\times\La^+$
	\end{itemize}
	We denote $\La^-\prec\La^+$ if there is a concordance from $\La^-$ to $\La^+$.
\end{defi}
	
Up to now it is not known if this relation induces a partial order on the set of Legendrian isotopy classes of Legendrian knots. It is of course well defined, transitive and reflexive but it is not known if $\La^-\prec\La^+$  and $\La^+\prec\La^-$  implies that $\La^-$ and $\La^+$  are Legendrian isotopic. There is a related notion of \emph{decomposable Lagrangian concordance } (denoted $\prec_{\operatorname{dec}}$) that are concordances built from elementary combinatorial moves. The projection to $\mathbb{R}$ of such a concordance only has critical points of index $0$ and $1$ hence are ribbon. It follows from recent work of Agol \cite{Ag} that if $\La^-\prec_{\operatorname{dec}}\La^+$  and $\La^+\prec_{\operatorname{dec}}\La^-$ then $\La^-$ and $\La^+$ are smoothly isotopic. Since Lagrangian concordances preserve both the Thurston-Bennequin number (TB) and the rotation number it implies that the relation $\prec_{\operatorname{dec}}$ is a partial order on the class of Legendrian knots whose isotopy classes are simple.
	
\begin{figure}[ht] 
	\labellist
	\pinlabel $m$ at 360 180
	\endlabellist
	\begin{center}\includegraphics[width=8cm]{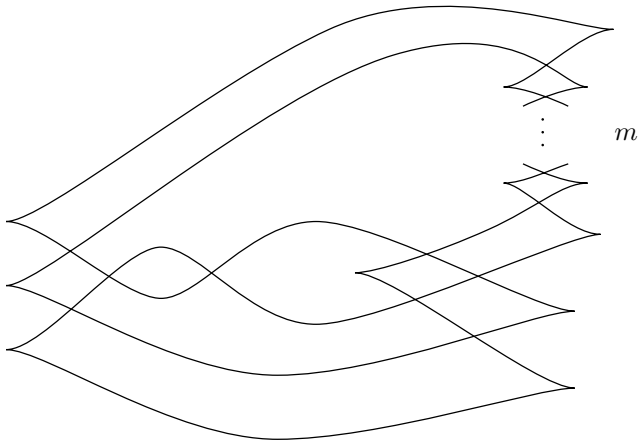}\end{center}
	\caption{Front projection of the Legendrian representative $\La_m$ of the pretzel knot $P(3,-3,-m)$.}
	\label{fig:pretzel}
\end{figure}

In this note we are concerned with Legendrian knots $\Lambda$ such that $\Lambda_0\prec \Lambda\prec \Lambda_0$ where $\Lambda_0$ is the maximal TB Legendrian unknot. The Legendrian knot $\Lambda_0$ is fillable by an exact Lagrangian disc. Thus it follows from the fact that a Lagrangian disc can be perturbed to a symplectic disc (because it is an open manifold) and a combination of results of Rudolph \cite{zbMATH03796828} and Boileau --Orevkov \cite{MR1836094} that the concordance realising  $\Lambda_0\prec \Lambda$ is ribbon. On the other hand in \cite[Theorem 3.2]{CNS} it is shown that the concordance realising $\Lambda\prec \Lambda_0$ cannot be decomposable, actually following the proof one observes that they prove that this concordance cannot be ribbon (which is also implied by a recent result of Zemke \cite{MR4024565} or the above mentioned result from \cite{Ag}).  In \cite{CNS} they give a list of Lagrangianly slice Legendrian knots for which the existence of a concordance to $\Lambda_0$ cannot be obstructed using Legendrian contact homology. The reason is that their Chekanov-Eliashberg algebras are stable tame isomorphic to the one of $\Lambda_0$ (see \cite[Appendix A]{EN} for explicit computation). The family in question consists of some $TB=-1$ Legendrian realisations of pretzel knots $P(3,-3,-m)$, for $m\geq4$. Denote $\La_m$ such a Legendrian knot, see Figure \ref{fig:pretzel} for its front projection.
There are concordances $\La_0\prec\La_m$ for all $m\geq4$, which can be easily constructed via some elementary moves on the front. 
	
In this note we make the following observation.
\begin{teo}\label{teo:concordance}
There is no concordance $\La_m\prec\Lambda_0$.
\end{teo}

Theorem \ref{teo:concordance} follows from a result of Eliashberg-Polterovitch \cite{EP} that we recall here (rephrased in the language that fits our purpose):
	
\begin{teo}\label{teo_elipo}
	Let $D$ be a filling of $\Lambda_0$, then $D$ is Hamiltonian isotopic to the standard filling of $\Lambda_0$.
\end{teo}

Indeed this result leads to obstructions for Legendrian knots to exist as middle slices of a concordance from $\Lambda_0$ to $\Lambda_0$. The symplectic flavour of the result leads to obstructions from Legendrian contact homology and variations of it, this has been used in \cite{Ch}, \cite{MR3762698} and \cite{W} for instance.  But the unknottedness of the disc in \cite{EP} has topological implications that to our knowledge have not been used in order to study Lagrangian cobordisms. Recall that a smooth knot $K$ is \emph{doubly slice} if there is an unknotted $2$-sphere $S$ in $\mathbb{R}^4$ such that $S\cap S^3=K$. As an immediate corollary of Theorem \ref{teo_elipo} we obtain:
\begin{teo}
	If there are concordances $\La_0\prec\La$ and $\La\prec\La_0$, then $\La$ is doubly slice.
\end{teo}
	
Now the main computation of this note shows that:
	
\begin{teo}\label{pretzel}
	The pretzel knots $P(3,-3,-m)$ for $m\geq4$ are not doubly slice.
\end{teo}
	
This implies Theorem \ref{teo:concordance}.

\begin{rem}
	For the case $m=3$, the knot $\La_3$ has topological type $m(9_{46})$ and it has been shown previously by the first author in \cite{Ch} that there is no concordance from $\La_3$ to $\La_0$. In this case computations using Legendrian contact homology or other modern tools for the case are still necessary because it turns out that the knot $m(9_{46})$ is doubly slice. The fact that the only one to be doubly slice in this family happens to be the only one to have a rich Legendrian contact homology algebra is both fortunate and puzzling.
\end{rem}

\begin{rem}
	In \cite{LM} the authors give an (almost complete) classification of doubly slice knots up to 12 crossings. From their work, the pretzel knots $P(3,-3,-4)=m(10_{140})$, $P(3,-3,-5)=11n_{139}$ and $P(3,-3,-6)=12n_{582}$ are not doubly slice.
\end{rem}
	
\section{Proof of Theorem \ref{pretzel}}
\subsection{Case $P(3,-3,-2k-1)$, $k\geq2$}
	
In this case the pretzel knots are odd because they have an odd number of parameters which are all odd. McDonald in \cite{McD}, using work of Issa and McCoy in \cite{MR4127867}, proves the following result:
	
\begin{teo}\cite[Theorem 3]{McD}
	For $K$ an odd pretzel knot, the following are equivalent:
\begin{itemize}
		\item $\Sigma_2(S^3,K)$ embeds in $S^4$,
		\item $K$ is a doubly slice pretzel knot,
		\item $K$ is a mutant of $P(a,-a,a,-a,\dots,a)$ for some odd $a$.
	\end{itemize}
\end{teo}
	
It is noted in \cite{McD} that it follows from the proof of  \cite[Theorem 1.11]{MR4127867} that an odd pretzel knot which is a mutant of $P(a,-a,a,-a,\dots,a)$ for some odd $a$ has for parameters a permutation of $(a,-a,a,-a,\dots,a)$. One concludes that the pretzel knots $P(3,-3,-2k-1)$ for $k\geq2$ are not doubly slice.

\subsection{Case $P(3,-3,-2k)$, $k\geq2$}
	
For these knots we obstruct doubly sliceness using the signature function as in \cite{LM} for $P(3,-3,-6)$. Given a knot $K$ and a Seifert matrix $A$ for $K$, the Levine-Tristram signature of $K$ is a function $\sigma_K:S^1\to\Z$ whose value at $\omega\in S^1$ is given by the signature of the matrix
\begin{alignat*}{1}
  (1-\omega)A+(1-\bar{\omega})A^T
\end{alignat*}
where the signature is the number of positive eigenvalues minus the number of negative eigenvalues. If $K$ is slice, then $\sigma_K(\omega)=0$ away of the roots of the Alexander polynomial of $K$. Moreover, if $K$ is doubly slice then $\sigma_K(\omega)$ vanishes for all $\omega\in S^1$. This is discussed in \cite[Proposition 2.2]{LM} and follows from the fact that when $K$ is doubly slice then there exists a Seifert matrix for $K$ that is hyperbolic (see \cite{MR290351}).
Seifert matrices for the knots $P(3,-3,-2k)$ have been computed in \cite[Section 3]{J}. We denote $A_k$ a Seifert matrix for $P(3,-3,-2k)$. We have:
\begin{alignat*}{1}
		A_k=\begin{pmatrix}
			-1&-1&0&0&-1&0\\
			0&-1&0&0&-1&0\\
			0&0&1&1&1&0\\
			0&0&0&1&1&0\\
			0&0&0&0&0&0\\
			0&0&0&0&1&-k
		\end{pmatrix}.
\end{alignat*}
It turns out that the Alexander polynomial of these knots is the same, given by $(t^2-t+1)^2$ and has for roots the sixth root of unity and its conjugate. One can then compute that $\sigma_K(\omega)\neq0$ for $K=P(3,-3,-2k)$ and $\omega=e^{i\frac{\pi}{3}}$. Indeed a direct computation shows that the matrix $(1-\omega)A_k+(1-\bar{\omega})A_k^T$ has rank $5$ and hence cannot have even signature (one can compute that its signature is actually $-1$).
Thus $K$ is not doubly slice which concludes both the proofs of Theorems \ref{pretzel} and \ref{teo:concordance}.\\

\paragraph{\textbf{Acknowledgement}}

This work was completed while the authors were at the Institute Mittag-Leffler for the program ``Frontiers of Quantitative Symplectic and Contact Geometry'', they want to thank the organisers and the staff from the institute for the stimulating environment they provided. The first author is partially supported by the ANR projects ENUMGEOM (ANR-18-CE40-0009) COSY (ANR-21-CE40-0002) and COSYDY (ANR-CE40-0014), the second author is supported by the grant KAW 2016.0198 from the Knut and Alice Wallenberg Foundation and the grant 2016-03338 from the Swedish Research Council.

	%
	%
	%
	%
	
\bibliographystyle{alpha}
\bibliography{ref.bib}	

\begin{thebibliography}{NRSS17}

\bibitem[Ago22]{Ag}
Ian Agol.
\newblock Ribbon concordance of knots is a partial order.
\newblock arXiv:2201.03626 [math.GT], 2022.

\bibitem[BO01]{MR1836094}
Michel Boileau and Stepan Orevkov.
\newblock Quasi-positivit\'{e} d'une courbe analytique dans une boule
  pseudo-convexe.
\newblock {\em C. R. Acad. Sci. Paris S\'{e}r. I Math.}, 332(9):825--830, 2001.

\bibitem[Cha15]{Ch}
Baptiste Chantraine.
\newblock Lagrangian concordance is not a symmetric relation.
\newblock {\em Quantum Topol.}, 6(3):451--474, 2015.

\bibitem[CNS16]{CNS}
Christopher Cornwell, Lenhard Ng, and Steven Sivek.
\newblock Obstructions to {L}agrangian concordance.
\newblock {\em Algebr. Geom. Topol.}, 16(2):797--824, 2016.

\bibitem[EN]{EN}
John~B. Etnyre and Lenhard Ng.
\newblock Legendrian contact homology in $\mathbb{R}^3$.
\newblock arXiv:1811.10966v3 [math.SG].

\bibitem[EP96]{EP}
Y.~Eliashberg and L.~Polterovich.
\newblock Local {L}agrangian {$2$}-knots are trivial.
\newblock {\em Ann. of Math. (2)}, 144(1):61--76, 1996.

\bibitem[IM20]{MR4127867}
Ahmad Issa and Duncan McCoy.
\newblock Smoothly embedding {S}eifert fibered spaces in {$S^4$}.
\newblock {\em Trans. Amer. Math. Soc.}, 373(7):4933--4974, 2020.

\bibitem[Jab10]{J}
Stanislav Jabuka.
\newblock Rational {W}itt classes of pretzel knots.
\newblock {\em Osaka J. Math.}, 47(4):977--1027, 2010.

\bibitem[LM15]{LM}
Charles Livingston and Jeffrey Meier.
\newblock Doubly slice knots with low crossing number.
\newblock {\em New York J. Math.}, 21:1007--1026, 2015.

\bibitem[McD20]{McD}
Clayton McDonald.
\newblock Doubly slice odd pretzel knots.
\newblock {\em Proc. Amer. Math. Soc.}, 148(12):5413--5420, 2020.

\bibitem[NRSS17]{MR3762698}
Lenhard Ng, Dan Rutherford, Vivek Shende, and Steven Sivek.
\newblock The cardinality of the augmentation category of a {L}egendrian link.
\newblock {\em Math. Res. Lett.}, 24(6):1845--1874, 2017.

\bibitem[Rud83]{zbMATH03796828}
Lee Rudolph.
\newblock Algebraic functions and closed braids.
\newblock {\em Topology}, 22:191--202, 1983.

\bibitem[Sum71]{MR290351}
D.~W. Sumners.
\newblock Invertible knot cobordisms.
\newblock {\em Comment. Math. Helv.}, 46:240--256, 1971.

\bibitem[Wu]{W}
Angela Wu.
\newblock Obstructing {L}agrangian concordance for closures of 3-braids.
\newblock arXiv:2201.08466v2 [math.SG].

\bibitem[Zem19]{MR4024565}
Ian Zemke.
\newblock Knot {F}loer homology obstructs ribbon concordance.
\newblock {\em Ann. of Math. (2)}, 190(3):931--947, 2019.

\end{thebibliography}
\end{document}